\documentclass[11pt]{article}
\usepackage{amsfonts}
\usepackage{mathrsfs}
\usepackage{amsmath}
\usepackage{amssymb}
\usepackage[mathscr]{eucal}
\usepackage{graphicx}
\def \vol{{\rm vol}}
\def \int{{\rm int}}
\def \la{\lambda}

\def \ov2{\overline}
\newtheorem{theorem}{\scshape \mdseries  Theorem}[section]

\newtheorem{corollary}[theorem]{\scshape \mdseries  Corollary}
\newtheorem{prop}[theorem]{\scshape \mdseries  Proposition}

\begin{document}
\title{\sf Spectral Condition for a Graph to be Hamiltonian with respect to Normalized Laplacian} \author{Yi-Zheng Fan$^{1,}$\thanks{Corresponding author.
Email: fanyz@ahu.edu.cn. Supported by National Natural Science
Foundation of China (11071002), Program for New Century Excellent
Talents in University, Key Project of Chinese Ministry of Education
(210091), Specialized Research Fund for the Doctoral Program of
Higher Education (20103401110002), Science and Technological Fund of
Anhui Province for Outstanding Youth (10040606Y33), Scientific
Research Fund for Fostering Distinguished Young Scholars of Anhui
University (KJJQ1001), Project for Academic Innovation Team of Anhui University (KJTD001B).},\ \ Gui-Dong Yu$^{1,2,}$\thanks{Email: yuguid@aqtc.edu.cn. Supported by NSF of Department of Education of Anhui Province (KJ2011A195) and Innovation Fund  for Graduates of Anhui University.}\\
  {\small  \it $1.$ School of Mathematical Sciences, Anhui University, Hefei 230601, P.R. China}\\
  {\small  \it $2.$ School of Mathematics \& Computation Sciences, Anqing Normal College, Anqing 246011, P.R. China}}
\date{}
\maketitle

\noindent {\bf Abstract:}
 Let $G$ be a graph and let $\Delta,\delta$ be the maximum and minimum degrees of $G$ respectively, where $\Delta/\delta<c<\sqrt{2}$ and $c$ is a constant.
 In this paper we establish a sufficient spectral condition for the graph $G$ to be Hamiltonian, that is, the nontrivial eigenvalues of the normalized Laplacian of $G$ are sufficiently close to $1$.


\noindent {\bf MR Subject Classifications:} 05C45, 05C50

\section{Introduction}

Let $G=(V,E)$ be a finite simple graph with vertex set $V=V(G)=\{v_1,v_2,\ldots,v_n\}$ and edge set $E=E(G)$.
The {\it adjacency matrix} of $G$ is defined to be a matrix $A=[a_{ij}]$ of order $n$, where $a_{ij}=1$ if $v_{i}$ is adjacent to $v_{j}$, and $a_{ij}=0$ otherwise.
Let $D$ be the diagonal matrix of order $n$ whose $(i,i)$-entry is $d_{v_{i}}$, the degree of the vertex $v_{i}$ of $G$.
The {\it signless Laplacian}, the {\it Laplacian}, and the {\it normalized Laplacian} of $G$ are respectively defined
  by $Q=D+A$, $L=D-A$ and $\mathcal{L}=D^{-1/2}LD^{-1/2}$ (for the last matrix we assume the graph contains no isolated vertices).

The graph $G$ is said to be {\it Hamiltonian} if there exists a cycle passing through all the vertices of $G$.
Such cycle is called a {\it Hamiltonian cycle} of $G$.
The question of deciding whether or not a given graph is Hamiltonian is a very difficult one;
  indeed, determining wether a given graph is Hamiltonian is NP-complete \cite{karp}.
Recently the spectral graph theory has been applied to this problem.
The sufficient spectral conditions are given for a graph having Hamiltonian paths
  or Hamiltonian cycles or being Hamilton-connected, in terms of spectral radius of a graph or its complement,
  with respect to the adjacency matrix or Laplacian or signless Laplacian; see Fiedler and Nikiforov \cite{fied}, Zhou \cite{zhou}, Yu and Fan \cite{yu}.
However, these conditions always imply the graph are very dense.

A breakthrough in studying Hamiltonicity occurred in 1975 when Koml\'os and Szemer\'edi \cite{kosz} showed that almost surely every random graph is Hamiltonian.
The technique involves the rotation of paths attributed to Posa \cite{posa}.
Krivelevich and Sudakov \cite{ks} established a sufficient condition for a $d$-regular graph to be Hamiltonian.
 They showed that if $\sigma$, the second largest absolute value of an eigenvalue of the adjacency matrix of a $d$-regular graph, satisfies
$$\sigma \le c\frac{(\log\log n)^{2}}{\log n(\log\log\log n)}d,\eqno(1.1)$$
for a constant $c$ and $n$ sufficiently large, then $G$ is Hamiltonian.
The condition (1.1) is not based on density conditions,
  rather it implies the graph is pseudo-random (the edge distribution resembles closely that of a truly random graph $G(n,d/n)$.

Using Laplacian of graphs, Butler and Chung \cite{but} established a sufficient condition for a graph $G$ being Hamiltonian.
They proved that if
$$|d-\mu_i|\leq c\frac{(\log\log n)^{2}}{\log n(\log\log\log n)}d,\eqno(1.2)$$
for $i \ne 0$, some constant $c$ and $n$ sufficiently large, then $G$ is Hamiltonian,
where $d$ is the average degree of $G$, and $0=\mu_0 \le \mu_1 \le \cdots \le \mu_{n-1}$ are the eigenvalues of the Laplacian of $G$.
The condition (1.2) implies the graph $G$ is almost regular, and in fact, pseudo-random.
If $G$ is regular, then (1.2) is exactly (1.1).

Mary Radcliffe \cite{mr} promoted the problem of finding sufficient conditions on the spectrum of the normalized Laplacian to ensure that a graph is Hamiltonian.
In this paper, we regard this problem and get the following result.
It can be seen the result also implies that of Krivelevich and Sudakov for regular graphs.

\begin{theorem} Let $G$ be a graph on $n$ vertices, $0=\la_0\leq \la_1 \le \cdots \la_{n-1}$ be the eigenvalues of the normalized Laplacian of $G$.
Assume that $\Delta/\delta<c<\sqrt{2}$ for some constant $c$, where $\Delta,\delta$ are
the maximum and minimum degrees of the vertices of $G$. If
$$|1-\la_i|\leq \frac{(\log\log n)^{2}}{7500\log n(\log\log\log n)},\eqno(1.3)$$
for $i\neq 0$ and $n$ sufficiently large, then $G$ is Hamiltonian.
\end{theorem}

{\bf Remark:}
We show two points on Theorem 1.1 by an example.
Let $G$ be the graph obtained from a complete graph $K_{n-1}$ on $n-1$ vertices by joining a new vertex with
$\beta:=\lceil \alpha (n-1) \rceil$ vertices of $K_{n-1}$, where $0 < \alpha <1$.
It is not too hard to show that:
$$
\max_{i \ne 0}|1-\la_i|=\frac{n-2+\beta+\sqrt{(n-2+\beta)^2+4(n-1)(n-2)(n-\beta-2)}}{2(n-1)(n-2)}\approx \frac{\sqrt{1-\alpha}}{\sqrt{n}} .
$$
So, this graph has the eigenvalues very tightly clustered near $1$ (i.e., even tighter than the bound in (1.3)).

(1) The constraint in Theorem 1.1 on the ratio of the maximal degree and minimal degree is necessary.
If taking $\alpha=\frac{1}{n-1}$, i.e., $G$ is $K_{n-1}$ with a pendant edge, surely $G$ is not Hamiltonian.
In this case $\frac{\Delta}{\delta}=n-1 \rightarrow \infty$.

(2) Theorem 1.1 applies more Hamiltonian graphs than Butler and Chung's result.
The condition (1.2) (or see Theorem 2.1 of \cite{but}) implies that $-\frac{1}{n}-\epsilon \le \frac{d_v}{d}-1 \le \epsilon$ for each vertex $v$, where $\epsilon=c\frac{(\log\log n)^{2}}{\log n(\log\log\log n)}$ .
So, when $n$ goes to infinity, $\frac{d_v}{d}-1 \rightarrow 0$, which implies the graph is almost regular.

For the above graph $G$, if taking $\alpha$ being a constant such that $\sqrt{2}/2 < \alpha < 1$,
then $\frac{\Delta}{\delta}<\sqrt{2}$,
Surely $G$ is Hamiltonian, which is consistent with our result.
However, $\left|\frac{\delta}{d}-1 \right| \rightarrow 1-\alpha \ne 0$.
So, using Butler and Chung's condition, we cannot decide whether it is Hamiltonian or not.

\section{Preliminaries}
Let $G$ be a graph, and let $X \subset V(G)$.
Denote by $\bar{X}$ be the complement of $X$ in $V(G)$, and by $N(X)$ the set of all vertices in $V \setminus X$ adjacent to some vertex in $X$.
The {\it volume} of $X$, denoted by $\vol(X)$, is defined as $\vol(X)=\sum_{v\in X}d_{v}$.
The volume of $G$ is denoted by $\vol(G)=\sum_{v\in G}d_{v}$.
For two subsets $X$ and $Y$ of $V$, we let $e(X,Y)$ be the number of edges with one endpoint in $X$ and one in $Y$, while $e(X)$ be the number of edges with both endpoints in $X$.

\begin{theorem} {\em \cite{fc}} \label{fan} Let $G$ be a graph on $n$ vertices, and let the eigenvalues $0=\la_0\leq \la_1 \le \cdots
\la_{n-1}$ of the normalized Laplacian of $G$ satisfy
$|1-\la_{i}|\leq\la$ for $i\neq 0$. Then for any two
subsets $X$ and $Y$ of the vertices in $G$,
$$\left|e(X,Y)-\frac{\vol(X)\vol(Y)}{\vol(G)}\right|\leq \la \frac{\sqrt{\vol(X)\vol(\bar{X})\vol(Y)\vol(\bar{Y})}}{\vol(G)}.$$
\end{theorem}

By Theorem \ref{fan}, we have the following conclusion immediately in terms of the maximum and minimum degrees.

\begin{corollary} \label{fandgr} Let $G$ be a graph on $n$ vertices with average degree $d$, and let the eigenvalues $0=\la_0\leq \la_1 \le \cdots
\la_{n-1}$ of the normalized Laplacian  of $G$ satisfy
$|1-\la_{i}|\leq\la$ for $i\neq 0$. Then for any two
subsets $X$ and $Y$ of the vertices in $G$,
\begin{align*}
e(X,Y) & \ge   \frac{\delta^{2}}{nd}|X||Y|- \frac{\la\Delta^{2}}{nd}\sqrt{|X|(n-|X|)|Y|(n-|Y|)}, \\
e(X,Y) & \le  \frac{\Delta^{2}}{nd}\left(|X||Y|+\la \sqrt{|X|(n-|X|)|Y|(n-|Y|)}\right).
\end{align*}
\end{corollary}

If we consider the case $X=\{v\}$ and $Y=V \setminus \{v\}$,
then Corollary \ref{fandgr} implies that
$$\frac{n-1}{n}\frac{\delta^{2}}{d}-\frac{\Delta^{2}}{d}\la\leq d_{v}\leq\frac{\Delta^{2}}{d}(1+\la).$$

\begin{corollary} \label{ex}
Let $G$ be a graph on $n$ vertices with average degree $d$,  and let the eigenvalues $0=\la_0\leq \la_1 \le \cdots \la_{n-1}$ of the normalized Laplacian of $G$ satisfy
$|1-\la_{i}|\leq\la$ for $i\neq 0$. Then for any subset $X$
of the vertices in $G$,
$$\frac{\delta^{2}}{2nd}|X|(|X|-1)-\frac{\la \Delta^{2}}{nd}|X|(n-|X|/2)\leq
e(X)\leq\frac{\Delta^{2}}{2nd}(|X|(|X|-1)+2\la|X|(n-|X|/2)).  $$
\end{corollary}

{\bf Proof:} Let $x=|X|$, $X'\subset X$ and $|X'|=\lfloor x/2\rfloor=x'$. Since
$$\sum_{X'\subseteq X \atop |X'|=x'}e(X',X\setminus X')=\binom{x-2}{x'-1}e(X,X),$$
 by the upper bound of $e(X,Y)$ in Corollary \ref{fandgr},  we have
\begin{align*}
\binom{x-2}{x'-1}e(X,X)& =\sum_{X'\subset X \atop |X'|=x'}e(X',X\setminus X')\\
& \leq \sum_{X'\subset X \atop |X'|=x'}\frac{\Delta^{2}}{nd}\left(|X'||X\setminus X'|+\la
\sqrt{|X'|(n-|X'|)|X\setminus X'|(n-|X\setminus X'|)}\right)\\
&=\binom{x}{x'}\frac{\Delta^{2}}{nd}\left(x'(x-x')+\la\sqrt{x'(n-x')(x-x')(n-x+x')}\right).
\end{align*}
So
\begin{align*}
e(X)&=\frac{1}{2}e(X,X)\\
& \leq  \binom{x}{x'}\binom{x-2}{x'-1}^{-1} \frac{\Delta^{2}}{2nd}\left(x'(x-x')+\la\sqrt{x'(n-x')(x-x')(n-x+x')}\right)\\
& \leq\frac{\Delta^{2}}{2nd}\left[x(x-1)+2\la x (n-\frac{x}{2})\right].
\end{align*}
 Similarly, by the lower bound of $e(X,Y)$ in Corollary \ref{fandgr},
we have
\begin{align*}
e(X)& =\frac{1}{2}e(X,X)\\
& \geq  \binom{x}{x'}\binom{x-2}{x'-1}^{-1} \left (\frac{\delta^{2}}{2nd}x'(x-x')-\frac{\la\Delta^{2}}{2nd}\sqrt{x'(n-x')(x-x')(n-x+x'))}\right)\\
& \geq\frac{\delta^{2}}{2nd}x(x-1)-\frac{\la \Delta^{2}}{nd}x(n-\frac{x}{2}).
\end{align*}   \hfill  $\blacksquare$

\begin{corollary}\label{basic}
Let $G$ be a graph on $n$ vertices with average degree $d$,  and let the
eigenvalues $0=\la_0\leq \la_1 \le \cdots \la_{n-1}$ of the
normalized Laplacian of $G$ satisfy
$|1-\la_i|\leq \la$ for $i\neq 0$. Further assume
that $\lambda<1/8$, $(\Delta/ \delta)^{2}\leq 2(n-1)/ n$, and that
$X,Y\subseteq V$. Then the following results hold:

{\em(a)} if $|X|<\la n$, then $e(X)\leq \frac{3\la \Delta^{2}}{2 d}|X|$;

{\em(b)} if $|X|<\la^{2} n$, then
$\displaystyle |N(X)|>\frac{(\frac{1}{2}-4\la)^{2}}{3\la^{2}}|X|$;

{\em(c)} if $|X|>\la^2 \Delta^4 n/\delta^4,$ then
$|N(X)|>\frac{n}{2}-|X|$;

{\em(d)} if $X\cap Y=\emptyset$ and $e(X,Y)=0$, then
$|X||Y|<\la^{2}\Delta^4 n^{2}/\delta^4$;

{\em(e)} $G$ is connected.
\end{corollary}

{\bf Proof:} For (a) we use Corollary \ref{ex} and the assumption to get
$$e(X)\leq\frac{\Delta^{2}}{2nd}(|X|(|X|-1)+2\la |X|(n-|X|/2))
\leq\frac{\Delta^{2}}{2nd}(\la n|X|+2\la n |X|)=\frac{3\la \Delta^{2}}{2 d}|X|.$$

For (b) if $|X|<\la^{2} n$,  then $|X|<\la n$. We use (a),
i.e., $e(X)\leq \frac{3\la \Delta^{2}}{2 d}|X|$, and
 the remark following Corollary \ref{fandgr},
 $$e(X,N(X))=\sum_{x\in X}d_{x}-2e(X)\geq \left(\frac{n-1}{n}\frac{\delta^{2}}{d}-\frac{\la\Delta^{2}}{d}\right)|X|
 -\frac{3\la\Delta^{2}}{d}|X|
=\left(\frac{n-1}{n}\frac{\delta^{2}}{d}-\frac{4\la\Delta^{2}}{d}\right)|X|. \eqno(2.1)$$

On the other hand, by Corollary \ref{fandgr},
\begin{align*}
e(X,N(X)) & \leq\frac{\Delta^{2}}{nd}\left(|X||N(X)|+\la \sqrt{|X|(n-|X|)|N(X)|(n-|N(X)|)}\right)\\
& \leq\frac{\Delta^{2}}{nd}|X||N(X)|+\frac{\la\Delta^{2}}{d} \sqrt{|X||N(X)|}.
\end{align*}
 If
$|N(X)|\leq\frac{(\frac{1}{2}-4\la)^{2}}{3\la^{2}}|X|$ then
we would have
\begin{align*}
\frac{\Delta^{2}}{nd}|X||N(X)|+\frac{\Delta^{2}}{d}\la \sqrt{|X||N(X)|}
& \leq\frac{\Delta^{2}}{nd}\frac{(\frac{1}{2}-4\la)^{2}}{3\la^{2}}|X|^{2}
 +\frac{\Delta^{2}}{d}\frac{\la(\frac{1}{2}-4\la)}{\sqrt{3}\la}|X|\\
&<\frac{\Delta^{2}}{d}\frac{(\frac{1}{2}-4\la)^{2}}{3}|X|+\frac{\Delta^{2}}{d}
\frac{\frac{1}{2}-4\la}{3/2}|X|\\
&\leq \frac{\Delta^{2}}{d}\left(\frac{1}{2}-4\la\right) |X|,
\end{align*}
using that $\la<1/8$, $(\Delta/ \delta)^{2}\leq 2(n-1)/ n$ in going to the last line, which is contradiction to (2.1), establishing (b).

For (c) letting $Y=V\setminus(X \cup N(X))$ and using Corollary \ref{fandgr}, we have
\begin{align*}
0=e(X,Y)& \geq\frac{\delta^{2}}{nd}|X||Y|-\la \frac{\Delta^{2}}{nd}\sqrt{|X|(n-|X|)|Y|(n-|Y|)}\\
& \geq\frac{\delta^{2}}{nd}|X||Y|-\la \frac{\Delta^{2}}{d}\sqrt{|X||Y|(1-\frac{|Y|}{n})},
\end{align*}
which upon rearranging gives
$$\frac{|Y|}{1-|Y|/n}\leq \frac{\la^{2}\Delta^4 n^{2}}{\delta^{4}|X|}<n.$$
This implies that $|Y|<n/2$ and hence $|N(X)|=n-|X|-|Y|>\frac n 2-|X|.$

For (d) again using Corollary \ref{fandgr}, we have
\begin{align*}
0=e(X,Y) & \geq\frac{\delta^{2}}{nd}|X||Y|-\la
\frac{\Delta^{2}}{nd}\sqrt{|X|(n-|X|)|Y|(n-|Y|)}\\
& >\frac{\delta^{2}}{nd}|X||Y|-\la \frac{\Delta^{2}}{d}\sqrt{|X||Y|};
\end{align*}
 and the result follows.

For (e), if $G$ is disconnected then $G$ has a connected component $X$ of size $|X|\leq n/2.$ Since $|N(X)|=\emptyset$, it follows from
part (c) that $|X|\leq\la^2 \Delta^4 n/\delta^4 \leq \frac{1}{8}\frac{4(n-1)^2}{n^2}\la n<\la n.$ We use (a),
i.e., $e(X)\leq \frac{3\la \Delta^{2}}{2 d} |X|$, and
 the remark following Corollary \ref{fandgr},
  \begin{align*}
  e(X,N(X))& =\sum_{x\in X}d_{x}-2e(X)\\
  & \geq \left(\frac{n-1}{n}\frac{\delta^{2}}{d}-\frac{\la\Delta^{2}}{d}\right)|X|
  -\frac{3\la\Delta^{2}}{d}|X|\\
  &=\left(\frac{n-1}{n}\frac{\delta^{2}}{d}-\frac{4\la\Delta^{2}}{d}\right)|X|\\
  & > \left(\frac{n-1}{n}\frac{\Delta^{2}n}{2(n-1)d}-\frac{\Delta^{2}}{2d}\right)|X|=0,
  \end{align*}
a contradiction. \hfill  $\blacksquare$

\section{Proof of Theorem 1.1}

The idea of the proof of Theorem 1.1 is to find a maximal path that
can be closed to create a cycle. Using the assumptions and Corollary
\ref{basic}, $G$ is connected, which implies that $G$ is Hamiltonian (if
not, there would be a vertex adjacent to some vertex in the cycle,
allowing us to create a path of longer length). The technique used
here is the rotation of the paths due to Posa \cite{posa}.

Let $P=(v_{1}, v_{2}, \ldots, v_{m})$ be a path of maximal length in
$G$. If $v_{m}$ is adjacent to $v_{i}$ (abbreviated $v_{i}\sim
v_{m}$) for some $i$, then another path of maximal length is given by
$P'=(v_{1}, \cdots,v_{i},v_{m},v_{m-1}, \cdots, v_{i+1})$. We say that
$P'$ is a {\it rotation} of $P$ with {\it fixed endpoint} $v_{1}$, {\it pivot}
$v_{i}$ and {\it broken edge} $v_{i}\sim v_{i+1}$. We can then rotate
$P'$ in a similar fashion to get a new path $P''$ of the same
length, and so on.

For $t\geq0$, let $S_{t}=\{$ $v\in V(P)\setminus \{v_{1}\}: v$ is
the endpoint of a path obtainable from $P$ by at most $t$ rotations
with fixed endpoint $v_{1}$, and all broken edges in $P\}$

\begin{prop}{\em\cite{ks}} \label{expan}For $t\geq0$,
$|S_{t+1}|\geq\frac{1}{2}|N(S_{t})|-\frac{3}{2}|S_{t}|$.
\end{prop}

Let
\begin{align*}
\la & =\frac{(\log\log n)^{2}}{7500\log n(\log\log\log n)};\\
t_{0}& =\left\lceil\frac{\log4\la^{2}n}{2\left(\log(1/(2\la)-4)-\log\sqrt{7}\right)}\right\rceil+2.
\end{align*}

By Corollary \ref{basic}(b), as long as $|S_{t}|<\la^{2} n$, then
$|N(S_{t})|>\frac{(\frac{1}{2}-4\la)^{2}}{3\la^{2}}|S_{t}|$,
and thus by Proposition \ref{expan},
$|S_{t+1}|>\frac{(\frac{1}{2}-4\la)^{2}}{6\la^{2}}|S_{t}|-\frac{3}{2}|S_{t}|$,
which implies
$$\frac{|S_{t+1}|}{|S_{t}|}>\frac{(\frac{1}{2}-4\la)^{2}}{7\la^{2}}.$$
In particular, using $\Delta/\delta < c<\sqrt{2}$, after at most
$t_{0}-2$ steps we have that
$|S_{t}|>\frac{\la^{2} n\Delta^{4}}{\delta^{4}}$.

By Corollary \ref{basic}(c) and Proposition \ref{expan} when taking one more step
we will have
$$|S_{t+1}|\geq\frac{1}{2}|N(S_{t})|-\frac{3}{2}|S_{t}|
\geq\frac{1}{2}\left(\frac{n}{2}-|S_{t}|\right)-\frac{3}{2}|S_{t}|\geq\frac{n}{4}-2|S_{t+1}|,$$
which implies $|S_{t+1}|\geq \frac{n}{12}$.

Let $Y=V\setminus (S_{t+1}\bigcup N(S_{t+1}))$, then $e(S_{t+1},Y)=0$.
By Corollary \ref{basic}(d), we have
$$|Y|<\frac{\la^{2}n^{2}\Delta^4/\delta^{4}}{|S_{t+1}|}
\le \frac{12\la^{2} n\Delta^{4}}{\delta^{4}}.$$
So,
$|N(S_{t+1})|=n-|S_{t+1}|-|Y|>n-|S_{t+1}|-\frac{12\la^{2} n\Delta^{4}}{\delta^{4}}.$

Again using Proposition \ref{expan}, we get
\begin{align*}
|S_{t+2}| & \geq\frac{1}{2}|N(S_{t+1})|-\frac{3}{2}|S_{t+1}| \\
& \geq\frac{1}{2}\left(n-|S_{t+1}|-\frac{12\la^{2} n\Delta^{4}}{\delta^{4}}\right)-\frac{3}{2}|S_{t+1}|\\
& = \frac{1}{2}(n-\frac{12\la^{2} n\Delta^{4}}{\delta^{4}})-2|S_{t+1}|\\
& > \frac{1}{2}n(1-48\la^2)-2|S_{t+2}|.
\end{align*}
So
$$|S_{t+2}|> \frac{1}{6}n(1-48\la^2) > \frac{n}{7}, \mbox{~i.e.~} |S_{t_{0}}|>\frac{n}{7}.$$

Let $B(v_{1})=S_{t_{0}}$ and $A_{0}=B(v_{1})\bigcup \{v_{1}\}.$ For
each $v\in B(v_{1})$ we can repeat the above argument to get $B(v)$,
$|B(v)|>n/7$, of endpoints of maximum length paths with endpoint
$v$. Note that each endpoint in $B(v)$ was obtained by at most
$2t_{0}$ rotations of $P$. So, for each $a\in A_{0}$,
$b\in B(a)$ there is a maximum length path $P(a,b)$ joining $a$ and
$b$ which is obtainable from $P$ by at most $\rho=2t_{0}$ rotations.

We return to the initial path $P$ and directed it. Since each
endpoint in $B(v_{1})$ is in $P$, we see $|P|\geq |B(v_{1})|>n/7$.
Then we can divide the path $P$ into $2\rho$ disjoint segments
$I_{1}, \cdots, I_{2\rho}$ each of length at least $\lfloor
n/14\rho\rfloor$. Since each path $P(a,b)$ is obtainable from $P$ by
at most $\rho$ rotations there are at least $\rho$ of the segments
untouched (but possibly traversed in the opposite direction). we
call each such segment {\it unbroken} in $P(a,b)$. These segments have an
absolute orientation induced by $P$, and another, relative to this
by $P(a,b)$ (where we direct that path from $a$ to $b$).

Let $$k=2\max\{1, \lceil 3000\rho\la\rceil\}.$$ We consider
sequence $\sigma=I_{i_{1}}, \cdots, I_{i_{k}}$ of $k$ unbroken
segments of $P$ which occur in this order in $P(a,b)$, where
$\sigma$ specifies not only the order of segments in $P(a,b)$ but
also their relative orientation. We say then that $P(a,b)$ contains
$\sigma$. Note that as $P(a,b)$  has at least $\binom{\rho}{k}$
sequences $\sigma$.

For a given $\sigma$ we denote by $L(\sigma)$ the set of all pairs
$a\in A_{0}$, $b\in B(a)$, for which the path $P(a,b)$ contains
$\sigma$.
The total number of possible sequences $\sigma$ is at most
$(2\rho)_{k}2^{k}$. Therefore by averaging we obtain that there
exists a sequence $\sigma_{0}$ for which
$$L(\sigma_{0})\geq\frac{n^{2}}{49}\frac{\binom{\rho}{k}}{(2\rho)_{k}2^{k}}
>\frac{n^{2}}{49}\left(\frac{\rho-k}{2\rho-k}\right)^{k}\frac{1}{k!2^{k}}.$$

It is easy to check that $k\leq \rho/2$ when $n$ sufficiently large.
Then $(\rho-k)/(2\rho-k)\geq1/3$, and it follows that there exists a
sequence $\sigma_{0}$ for which $|L(\sigma_{0})|\geq
n^{2}/(49k!6^{k})$. We fix such a sequence and denote
$$\alpha=\frac{1}{49k!6^{k}}.$$

Let $\hat{A}=\{a\in A_{0}:L(\sigma_{0})  \mbox{~contains at least~}
\alpha n/2 \mbox{~pairs with~} a \mbox{~as the first element}\}$.
Then $|\hat{A}|\geq \alpha n/2$. For each $a\in \hat{A}$, let
$\hat{B}(a)=\{b\in B(a): (a,b)\in L(\sigma_{0})\}.$ The definition
of $\hat{A}$ guarantees that $|\hat{B}|\geq \alpha n/2$¡£

Let $C_{1}$ be the union of the first $k/2$ segments of
$\sigma_{0}$, in the fixed order and with the fixed relative
orientation in which they occur along any of the paths $P(a,b)$,
$(a,b)\in L(\sigma_{0})$. Let $C_{2}$ be the union of the last $k/2$
segments of $\sigma_{0}$. Note that for $i=1,2$,
$$|C_{i}|\geq\frac{k}{2}\left\lfloor
\frac{n}{14\rho}\right\rfloor
\geq 3000\rho\la\left\lfloor
\frac{n}{14\rho}\right\rfloor>200n\la.  \eqno(3.1)$$

Given a path $P$ and a set $S\subset V(P)$, a vertex $v\in S$ is
called an {\it interior point} of $S$ with respect to $P$ if both
neighbors of $v$ along $P$ are in $S$. The set of all interior
points of $S$ will be denoted by $\int(S)$.

\begin{prop}
The set $C_{1}$ contains a subset $C'_{1}$ with
$|\int(C'_{1})|\geq nk/(60\rho)$ so that every vertex $v\in
C'_{1}$ has at least $48\la d$ neighbors in $\int(C'_1)$.
A similar statement holds for $C_{2}$.
\end{prop}

{\bf Proof:} We start with $C'_1=C_{1}$ and as long as there
exists a vertex $v_{j}\in C'_1$ for which has less than $48\la
d$ neighbors in $\int(C'_1)$, we delete $v_{j}$ and repeat. If
this procedure continued for $r=|C_{1}|/8$ steps then we get a
subset $R=\{v_{1}, v_{2}, \cdots, v_{r}\}$, so that
$$|\int(C'_1)|\geq
|\int(C_{1})|-3r=(1-o(1))|C_{1}|-3r
>|C_{1}|/2 \geq nk/(60\rho)$$ and
$$e(R, \int(C'_1))\leq 48\la dr=6\la d|C_{1}|.  \eqno(3.2)$$

But according to Corollary \ref{fandgr} and (3.1),
\begin{align*}
e(R, \int(C'_1)) & \geq
\frac{\delta^{2}}{nd}|R||\int(C'_1)|-
\frac{\la\Delta^{2}}{nd}\sqrt{|R|(n-|R|)|\int(C'_1)|(n-|\int(C'_1)|)}\\
& \geq
\frac{\delta^{2}}{nd}|R||\int(C'_1)|-
\frac{\la\Delta^{2}}{d}\sqrt{|R||\int(C'_1)|}\\
& \geq
\frac{\delta^{2}}{nd}\frac{|C_{1}|^2}{16}-
\frac{\la\Delta^{2}}{d}\sqrt{\frac{|C_{1}|^2}{16}}
>\frac{\delta^{2}}{nd}\frac{200n\la|C_{1}|}{16}-\frac{\la\Delta^{2}|C_1|}{4d}\\
&=\la d|C_{1}|\left(\frac{25\delta^2}{2d^2}-\frac{\Delta^2}{4d^2}\right)
\ge \la d|C_{1}|\left(\frac{25\delta^2}{2\Delta^2}-\frac 1 4 \right)\\
& >6\la d|C_{1}|.
\end{align*}
using that $\delta^2/\Delta^2>1/2$ in going to the last line,
which is contradiction to (3.2). So, the result follows. \hfill $\blacksquare$

We fix the obtained sets $C'_1$ and $C'_2$.

\begin{prop} There is a vertex $\hat{a}\in \hat{A}$ connected by an edge to $\int(C'_1)$.
Similarly there is a vertex $\hat{b}\in \hat{B}(\hat{a})$ connected
by an edge to $\int(C'_2)$.
\end{prop}

{\bf Proof:} Recall that $|\hat{A}|\geq\frac{\alpha n}{2}$, and
$|\int(C'_1)|\geq nk/(60\rho)$. Therefore, by Corollary \ref{basic}(d),
the claim will follows if we will show that
$\frac{\alpha n}{2} \frac{nk}{60\rho} \gg \frac{\Delta^4 \la^{2} n^{2}}{\delta^4}$, or
 (substituting the value of $\alpha$) $\frac{\delta^4}{\Delta^4\la^{2}\rho}\gg
 5880(k-1)!6^{k}$.

Consider first the case $3000\rho\la\geq1$. In this case,
\begin{align*}
k& =2(1+o(1))3000\la\rho=6000(1+o(1))\frac{\la\log \la^2 n}{\log(1/\la)}\\
& \leq 6000(1+o(1))\left(-2\la+\frac{\log n}{(1/\la) \log(1/\la)}\right)\\
& =6000(1+o(1))\frac{\log n}{\frac{7500\log n(\log\log\log n)}{(\log\log n)^{2}}\log\log n}\\
&=\frac{0.8(1+o(1))\log\log n}{\log\log\log n},
\end{align*}
and thus $5880(k-1)!6^{k}<(\log n)^{0.9}$.
On the other hand, as $\delta/\Delta> 1/\sqrt{2}$,
\begin{align*}
\frac{\delta^4}{\Delta^4\la^{2}\rho}
 & \geq \frac{1}{4\la^{2}}\frac{\log (1/\la)}{(1+o(1))\log \la^2 n}\\
 &\geq \frac{\log^2 n(\log\log\log n)^2}{(\log\log n)^4}\frac{\log \log n}{(1+o(1))\log n}\\
 &>\frac{(1+o(1))\log n}{(\log \log n)^{3}}\\
 &\gg(\log n)^{0.9},
 \end{align*} as required.

In the second case, $3000\rho\la<1$, we get $k=2$, then the
expression $(k-1)!6^{k}$ is an absolute constant, while
$\frac{\delta^4}{\Delta^4\la^{2}\rho} \geq \frac 1 4 \frac{1}{\rho\la} \frac{1}{\la}
\ge \frac{750}{\la} \rightarrow\infty$.
 The Proposition follows.  \hfill  $\blacksquare$

Now, let $x$ be a vertex separating $C'_1$ and $C'_2$
along $P(\hat{a},\hat{b})$, we consider two half path $P_{1}$ and
$P_{2}$ obtained by splitting $P(\hat{a},\hat{b})$ at $x$.
Consider $P_{1}$ firstly. Let $T_{i}=\{ v\in C'_1\setminus \{x\}: v$ is
the endpoint of a path obtainable from $P_{1}$ by $i$ rotations with
fixed endpoint $x$, all pivots in $\int(C'_1)$ and all broken
edges in $P_{1}\}$.

\begin{prop} There exists an $i$ for which $|T_{i}|\geq\la n(\Delta/\delta)^{2}$.
\end{prop}

{\bf Proof:} It is enough to prove that there exits a sequence of
sets $U_{i}\subseteq T_{i}$ such that $U_{1}=1$ and
$|U_{i+1}|=2|U_{i}|$, as long as $|U_{i}|<\la
n(\Delta/\delta)^{2}$. Note that according to Proposition 3.3
$\hat{a}$ has a neighbor in $\int(C'_1)$, and therefore
$T_{1}\neq\emptyset$. Note also that if we perform a rotation a
vertex from $\int(C'_1)$ and broken edge in $P_{1}$, then the
resulting endpoint is in $C'_1$.

Suppose we have found sets $U_{1}, \cdots, U_{i}$ as state above,
and still $|U_{i}|<\la n(\Delta/\delta)^{2}$.
We first show that
$$|T_{i+1}|\geq\frac{1}{2}|N(U_{i})\cap \int(C'_1)|-\frac{3}{2}\sum_{j=1}^{i}|U_{j}|.$$
Let
$$T=\{k \ge 1: v_k \in N(U_i) \cap \int(C'_1), v_{k-1},v_k,v_{k+1} \notin \cup_{j=1}^i U_j\}.$$
Consider a vertex $v_k$ with $k \in T$.
Then $v_k$ has a neighbor $w \in U_i$ which is also a interior vertex of $C'_1$.
So there exists a path $Q$ with $w$ as an endpoint, obtained from $P_{1}$ by $i$ rotations with
fixed endpoint $x$.
As $ v_{k-1},v_k,v_{k+1} \notin \cup_{j=1}^i U_j$, both edges $(v_{k-1},v_k)$ and $(v_k,v_{k+1})$ are still present in $Q$.
Rotating $Q$ with a pivot $v_k$ and one of the edges $(v_{k-1},v_k)$ and $(v_k,v_{k+1})$ as
a broken edge will put one of $v_{k-1},v_{k+1}$, say $v_{k-1}$ in $T_{i+1}$.
The only other vertex that possible cause $v_{k-1}$ to be put into $T_{i+1}$ is $v_{k-2}$ if
$k-2 \in T$.
Therefore,
$$ |T_{i+1}| \ge \frac 1 2 |T| \ge \frac 1 2 (|N(U_{i})\cap \int(C'_1)|-3|\sum_{j=1}^{i}U_{j}|)
\ge \frac{1}{2}|N(U_{i})\cap \int(C'_1)|-\frac{3}{2}\sum_{j=1}^{i}|U_{j}|.$$

As $\sum_{j=1}^{i}|U_{j}|<2|U_{i}|$, the claim will follow if we
prove that $|N(U_{i})\cap \int(C'_1)|\geq10|U_{i}|.$
Since $U_{i}\subset C'_1$, every vertex $u\in U_{i}$ has at
least $48\la d$ neighbors in $\int(C'_1)$. Therefore $e(U_{i}, \int(C'_1))\geq 48\la d|U_{i}|$.
Let
$W_{i}=N(U_{i})\bigcap \int(C'_1)$. If
$|W_{i}|<10|U_{i}|$, then by Corollary \ref{fandgr} we have
\begin{align*}
e(U_{i}, W_{i})& \leq\frac{\Delta^{2}}{nd}\left(|U_{i}||W_{i}|+\la
\sqrt{|U_{i}|(n-|U_{i}|)|W_{i}|(n-|W_{i}|)}\right)\\
& \leq\frac{\Delta^{2}}{nd}|U_{i}||W_{i}|+\frac{\la\Delta^{2}}{d}
\sqrt{|U_{i}||W_{i}|}<\frac{10\Delta^2 |U_i|^2}{nd}+\frac{\sqrt{10}\la \Delta^2|U_i|}{d}\\
&=\frac{\Delta^2}{d^2}\la d|U_{i}|\left(\frac{10|U_i|}{\la n}+\sqrt{10}\right)
< \frac{\Delta^2}{\delta^2}\la d|U_{i}|\left(\frac{10\la n(\Delta/\delta)^{2}}{\la n}+\sqrt{10}\right)\\
& < 2 \la d |U_i| (20+\sqrt{10})<48\la d |U_i|,
\end{align*}
a contradiction. Therefore $|W_{i}| \ge 10|U_{i}|$, as desired.  \hfill  $\blacksquare$

Hence, the set $V_{1}$ of endpoints of all rotations of $P_{1}$ has
cardinality $|V_{1}|\geq\la n(\Delta/\delta)^{2}$. Similarly
the set $V_{2}$ of endpoints of all rotations of $P_{2}$ also has
cardinality $|V_{2}|\geq\la n(\Delta/\delta)^{2}$. Then,
$|V_{1}||V_{2}|\geq(\Delta/\delta)^{4}\la^{2} n^{2}$, by
Corollary \ref{basic}(d) there is an edge connecting $V_{1}$ and $V_{2}$ and
thus closing the cycle. As $G$ is connected by Corollary \ref{basic}(e),
this cycle is a Hamilton cycle. This completes the proof of Theorem
1.1. \hfill  $\blacksquare$

\end{document}